\documentclass{amsart}
\usepackage{amssymb,amsmath,amsthm,xspace}
\usepackage[small]{diagrams}
\usepackage{ifpdf}
\ifpdf\usepackage[pdftex]{hyperref}%
\else\usepackage{hyperref} \fi

\def\mat#1{\ensuremath{#1}\xspace}

\def\makemath#1#2#3{%#1-new command, #2-new box, #3-math expression
\newsavebox{#2}
\sbox{#2}{\ensuremath{#3}}
\def#1{\usebox{#2}\xspace}}

\def\cP{\mat{\mathbb{P}}}   %projective space
\def\cZ{\mat{\mathbb{Z}}}   %ring of integers

\def\lB{\mat{\mathcal{B}}}

\def\lE{\mat{\mathcal{E}}}

\def\lH{\mat{\mathcal{H}}}

\def\lO{\mat{\mathcal{O}}}

\def\lS{\mat{\mathcal{S}}}

\makemath{\Psi}{\Psibox}{\Psi}
\makemath{\Phi}{\Phibox}{\Phi}

\def\la{\mat{\lambda}}

\def\hi{\mat{\chi}}

\def\mrm@#1{\mat{\mathrm{#1}}}

\def\DMO{\DeclareMathOperator}

\DMO{\Hom}{Hom}
\DMO{\lHom}{\lH\mathit{om}}
\DMO{\Ext}{Ext}
\DMO{\lExt}{\lE\mathit{xt}}
\DMO{\End}{End}
\DMO{\Aut}{Aut}
\DMO{\Fun}{Fun}
\DMO{\Tor}{Tor}
\DMO{\ext}{ext}

\DMO{\Ob}{Ob}
\DMO{\Mor}{Mor}
\DMO{\im}{im}
\DMO{\coim}{coim}
\DMO{\coker}{coker}
\DMO{\Arr}{Arr}

\DMO{\Id}{Id}
\DMO{\add}{add} % splitting of idempotents (karoubinization)
\DMO{\ind}{ind} % category of ind-objects
\DMO{\pro}{pro} % category of pro-objects
\DMO{\Map}{Map}
\DMO{\Iso}{Iso}
\DMO{\Isom}{Isom}

\DMO{\Presh}{Presh}

\DMO\coalg{Coalg}
\DMO{\Rep}{Rep}

\DMO{\Mod}{Mod}
\DMO{\rad}{rad}
\DMO{\soc}{soc}
\DMO{\ann}{ann}

\DMO{\Spec}{Spec}
\DMO{\spec}{Spec}
\DMO{\Proj}{Proj}
\DMO{\supp}{supp}
\DMO{\Coh}{Coh}
\DMO{\coh}{Coh}
\DMO{\Qcoh}{QCoh}
\DMO{\QCoh}{QCoh}
\DMO{\Pic}{Pic}
\DMO{\Div}{Div}
\DMO{\ch}{ch}
\DMO{\Hilb}{Hilb}
\DMO{\Fitt}{Fitt}
\DMO{\Quot}{Quot}

\DMO{\Gras}{Gr}
\DMO{\Flag}{Flag}

\DMO{\cone}{cone}
\DMO{\Tw}{Tw}
\DMO{\rank}{rk}
\DMO{\rk}{rk}
\DMO{\codim}{codim}
\DMO{\cov}{cov}
\DMO{\sgn}{sgn}
\DMO{\td}{td}
\DMO{\GL}{GL}
\DMO{\SL}{SL}

\DMO\Der{Der}
\DMO\der{Der}
\DMO\coder{Coder}
\DMO{\diag}{diag}
\DMO{\HMod}{HMod} %the homotopy category of modules over DGC
\DMO{\ad}{ad}
\DMO*{\colim}{colim}
\DMO*{\hocolim}{hocolim}
\DMO*{\holim}{holim}
\DMO{\Ho}{Ho}
\DMO{\har}{char}
\DMO{\sk}{sk}
\DMO{\cosk}{cosk}
\DMO{\Gal}{Gal}
\DMO{\tr}{tr}
\DMO{\Tr}{Tr}
\DMO{\Sh}{Sh}
\DMO{\Is}{Is} %Isometries
\DMO{\Hol}{Hol} %Holomorphic automorphisms
\DMO{\Lie}{Lie} %Lie algebra of a group
\DMO{\Res}{Res} %restriction
\DMO{\irr}{irr} %
\DMO{\Irr}{Irr} %
\DMO{\Exp}{Exp} %
\DMO{\Log}{Log} %
\DMO{\mult}{mult} %
\DMO{\height}{ht}

\def\iso{\simeq}
\def\ts{\otimes}

\def\wtl#1{\widetilde{#1}}

\def\mb#1{\mat{\mathbf{#1}}}
\def\br{\linebreak}
\def\arrowsFromDiagrams{
\newarrow{ShortTo}{}{}-->
\def\artest{{\!\!\!\!\!\!\!\rShortTo\!}}
\def\ar{{\:\vrule depth -.52ex height .60ex width 0.85em\;\!\!\rhla\,}}
\def\emb{{\:\rthooka\!\!\!\ar}}
\def\mto{{\:\vrule height .9ex depth -.2ex width .04em\!\!\!\;\ar}}
\def\arr{\rTo} % to use with indexes
\def\embb{\rInto} % to use with indexes
\newarrow{Eq}=====
\def\rrarr{\pile{\rTo\\ \rTo}}
\def\lrarr{\pile{\rTo\\ \lTo}}   }

\newif\ifukr\ukrfalse
\newif\ifrus\rusfalse

\def\theorems{
\newtheorem{prp}{\ifukr Пропозиція \else \ifrus Предложение \else Proposition\fi\fi}[section]
\newtheorem{proposition}[prp]{\ifukr Пропозиція \else \ifrus Предложение \else Proposition\fi\fi}

\newtheorem{theorem}[prp]{\ifukr Теорема \else \ifrus Теорема \else Theorem\fi\fi}

\newtheorem{lemma}[prp]{\ifukr Лема \else \ifrus Лемма \else Lemma\fi\fi}

\newtheorem{corollary}[prp]{\ifukr Висновок \else \ifrus Следствие \else Corollary\fi\fi}

\theoremstyle{definition}

\newtheorem{definition}[prp]{\ifukr Означення \else \ifrus Определение \else Definition\fi\fi}

\newtheorem{remark}[prp]{\ifukr Зауваження \else \ifrus Замечание \else Remark\fi\fi}
}

\def\line#1#2#3{\mat{#1_{#2},\dots,#1_{#3}}}
\theorems\arrowsFromDiagrams
\def\mch{M_{\mathrm{ch}}}
\def\mcyc{M_{\mathrm{cyc}}}

\begin{document}
\title[]{Classification of semistable sheaves on a rational curve with one node}%
\author{Sergey Mozgovoy}%
\address{Institut f\"ur Mathematik, Johannes Gutenberg-Universit\"at Mainz,
55099 Mainz, Germany.}%
\email{mozgov@mathematik.uni-mainz.de}%

%\thanks{}%
%\keywords{}%
%\dedicatory{}%
%\commby{}%
\maketitle
\begin{abstract}
We classify (semi)stable sheaves on a rational curve with one
node. The results are based on the classification of
indecomposable torsion-free sheaves due to Drozd and Greuel
\cite{DG}, where the sheaves are described in terms of certain
combinatorial data. We translate the condition of (semi)stability
into this combinatorial language and solve the so obtained
problem.
\end{abstract}

\section{Introduction}
A rational curve with one node is among those few examples of
singular curves, where a complete classification of indecomposable
vector bundles is possible. The classification was done by Drozd
and Greuel in~\cite{DG} in terms of certain combinatorial objects,
which are easy to handle. Using their technique one can also
classify all indecomposable torsion-free sheaves.

A next natural question would be to describe all semistable and
stable sheaves.
%However, until now it has been done only in the following particular cases.
However, until now only some partial results in this direction
were known. The problem was solved for the sheaves of degree $0$
by Burban and Kreu\ss ler \cite{BK1,BK2} by reducing it to the
classification of torsion sheaves made by Gelfand and Ponomarev
\cite{GP1}. In the case of coprime degree and rank Burban
\cite{Bu} classified stable locally free sheaves by detecting
those sheaves that have a one-dimensional endomorphism ring.

In this paper we give a classification of indecomposable
(semi)stable sheaves (of nonzero rank) on a rational curve with
one node over an algebraically closed field of characteristic $0$.
In order to do this we introduce and analyze certain combinatorial
objects --- chains and cycles (see Definition \ref{def:ch:cyc}),
which are used for the classification of indecomposable
torsion-free sheaves. With any aperiodic cycle \mb a one
associates an indecomposable locally free sheaf $\lB(\mb a)$ (see
\cite{DG,BK1} and Section \ref{sec:ss:sh}) and with any chain \mb
b one associates an indecomposable non-locally free sheaf $\lS(\mb
b)$ (see \cite{BK1} and Section \ref{sec:ss:sh}). The conditions
of (semi)stability of $\lB(\mb a)$ and $\lS(\mb b)$ imply certain
conditions on the cycle \mb a and the chain \mb b. We will call
these conditions the conditions of (semi)stability of cycles and
chains, respectively (see Section \ref{sec:ss:ch:cyc} for precise
definitions). One of the main results of the paper is

\begin{theorem}\label{intro:main:thr1}
Given an aperiodic cycle \mb a, the sheaf  $\lB(\mb a)$ is
(semi)stable if and only if the cycle $\mb a$ is (semi)stable.
Given a chain $\mb b=(\line b1r)$, the sheaf $\lS(\mb b)$ is
(semi)stable if and only if the chain
$(b_1+1,b_2,\dots,b_{r-1},b_r+1)$ is (semi)stable.
\end{theorem}

This means that we only need to classify the (semi)stable chains
and cycles. This is a purely combinatorial problem, and it has the
following solution. Let $\mcyc^{ss}(r,d)$ (respectively,
$\mcyc^s(r,d)$, $\mch^{ss}(r,d)$, $\mch^s(r,d)$) be the set of all
aperiodic semistable cycles (respectively, aperiodic stable
cycles, semistable chains, stable chains) of rank $r$ and degree
$d$.

\begin{theorem}\label{intro:main:thr2}
Let $r\in\cZ_{>0}$ and $d\in\cZ$. Then
\begin{enumerate}
    \item There is a natural bijection between $\mch^{ss}(r,d)$ and
$\mch^{ss}(r,d+r)$ and if $0<d<r$ then there is a natural
bijection between $\mch^{ss}(r,d)$ and $\mch^{ss}(d,d-r)$. As a
corollary, there is a bijection between $\mch^{ss}(r,d)$ and
$\mch^{ss}(h,0)$, where $h=\gcd(r,d)$. The same assertions hold
for stable chains, stable cycles and semistable cycles.
    \item The set $\mch^{ss}(r,d)$
is finite and non-empty. If $r$ and $d$ are coprime then
$\mch^{s}(r,d)=\mch^{ss}(r,d)$ and it contains just one element.
Otherwise, $\mch^s(r,d)$ is empty. The same assertions hold for
cycles.
\end{enumerate}
\end{theorem}

Thus the description of $\mch^{ss}(r,d)$ is reduced to the
description of $\mch^{ss}(h,0)$, $h=\gcd(r,d)$, and the latter is
given in Proposition~\ref{prp:sschains:r0} (for the analogous
classification of cycles, see Proposition~\ref{prp:sscycles:r0}).
Among other things, we prove that for $h>1$ there are no stable
chains (cycles) and for $h=1$ there is just one semistable chain
(cycle) which is actually stable. This implies the second part of
Theorem \ref{intro:main:thr2}.

In Section \ref{sec:ss:sh} we recall the definition of
indecomposable sheaves associated to chains and cycles. Writing
down the (semi)stability condition for these sheaves we get
certain conditions on chains and cycles, which we call the
(semi)stability conditions on chains and cycles.

In Section \ref{sec:ss:ch:cyc} we describe basic properties of
chains and cycles. We define (semi)\br stable chains and cycles,
analyze their structure and give basic reduction methods.
Altogether this allows us to classify the (semi)stable chains and
cycles.

In Section \ref{sec:ss:sh:classif} we use the classification from
Section \ref{sec:ss:ch:cyc} to prove that for any (semi)stable
chain or cycle, the associated sheaf is also (semi)stable.
Together with the results from Section \ref{sec:ss:sh} this proves
that the conditions of (semi)stability of chains and cycles are
necessary and sufficient for the (semi)stability of the
corresponding sheaves. As a corollary, we prove in particular that
any indecomposable semistable sheaf is homogeneous, i.e., all
stable factors of its Jordan-H\"older filtration are isomorphic.
This was proved in \cite{FMW} for sheaves of degree $0$.

My cordial thanks go to N.Sidorova. Proposition
\ref{lmm:ch:reduction} is a result of our communication. I would
like to thank also I.Burban, M.Lehn, C. Pillau, C. Sorger, I.Yudin
for many useful comments.

\section{Semistable sheaves}
\label{sec:ss:sh}
Let $C$ be a rational curve with one node over an algebraically
closed field $k$ of characteristic $0$. Let $\pi:\wtl C\ar C$ be
its normalization ($\wtl C\iso\cP^1$). Given a torsion-free sheaf
$F$ over $C$ of rank $r$ and degree $d$ ($\deg
F=\hi(F)-r\hi(\lO_C)=\hi(F)$), we will say that $F$ is of type
$(r,d)$. We will denote the set of torsion-free, indecomposable
sheaves of type $(r,d)$ by $\lE(r,d)$. Such sheaves were
classified by Drozd and Greuel \cite{DG} in terms of certain
combinatorial data. Our aim is to express the conditions of
stability and semistability of sheaves in $\lE(r,d)$ in the
language of this combinatorial data. As a result we will get a
classification of stable and semistable sheaves.

Let us give a description of sheaves in $\lE(r,d)$ according to
\cite{DG}. Let $p,\ p^*$ be preimages of a singular point in $C$
under $\pi$. For any line bundle $L\iso \lO(n)$ over $\wtl C$ we
fix once and for all the bases of the fibers $L(p)$ and $L(p^*)$.
To make possibly few choices we do this in the following way. Fix
some section $s$ of $\lO(1)$ having zero in some point different
from $p$ and $p^*$. Then $s^n$ will induce nonzero elements of the
fibers of $\lO(n)$ over $p$ and $p^*$, giving the necessary bases.

\begin{definition}\label{def:ch:cyc}
Define a chain to be a finite sequence of integers. Define a cycle
to be an equivalence class of chains, where the equivalence is
generated by relations
$$(a_1,a_2,\dots,a_r)\sim (a_2,a_3,\dots,a_r,a_1).$$
We will usually write representing sequences instead of the
corresponding cycles.
\end{definition}

Given a finite sequence of integers $\mb a=(\line a1r)$, a natural
number $m$ and an element $\la\in k^*$, we construct the vector
bundle $\lB(\mb a,m,\la)$ over $C$ in the following way (see
\cite{DG} or \cite{BK1} for more formal description). Consider the
sheaves $B_i=\lO(a_i)^{\oplus m}$ over $\wtl C$, then take the
direct image $\pi_*$ of their sum and make the following
identifications over the singular point: glue $B_1(p^*)$ with
$B_2(p)$, glue $B_2(p^*)$ with $B_3(p)$ and so on up to
identification of $B_r(p^*)$ with $B_1(p)$. The gluing matrices
(with respect to the above chosen bases) are defined to be unit
matrices except the matrix gluing $B_r(p^*)$ with $B_1(p)$ which
is defined to be a Jordan block of size $m$ with an eigenvalue
\la.

\begin{remark}
Note that if $\mb a=(\line a1r)$ and $\mb {a'}=(\line a2r,a_1)$ is
its cyclic shift, then $\lB(\mb a,m,\la)\iso\lB(\mb{a'},m,\la)$.
This means that the vector bundle $\lB(\mb a,m,\la)$ is determined
by the cycle $\mb a$ (together with $m\in\cZ_{>0}$, $\la\in k^*$).
\end{remark}

\begin{theorem}[see \cite{DG}]
The sheaves $\lB(\mb a,m,\la)$ with aperiodic (see Definition
\ref{def:aperiodic}) cycles \mb a describe all indecomposable
locally free sheaves over $C$. Different cycles induce
non-isomorphic sheaves.
\end{theorem}

In particular, consider the cycle $\mb 0=(0)$ and define
$F_m:=\lB(\mb 0,m,1)$. One can show that $F_1\iso \lO_C$ and there
is an exact sequence
$$0\ar \lO_C\ar F_m\ar F_{m-1}\ar 0.$$

\begin{lemma}[see \cite{Yu}]
There are isomorphisms
$$\lB(\mb a,m,\la)\iso \lB(\mb a,1,\la)\ts F_m, \qquad
\lB(\mb a,1,\la^r)\iso \lB(\mb a,1,1)\ts \lB(\mb 0,1,\la),$$ where
$r$ is the length of \mb a.
\end{lemma}

We will denote $\lB(\mb a,1,1)$ by $\lB(\mb a)$.

\begin{corollary}\label{crl:lf:sheaves:reduction}
The sheaf $\lB(\mb a,m,\la)$ is semistable if and only if $\lB(\mb
a)$ is semistable. The sheaf $\lB(\mb a,m,\la)$ is stable if and
only if $m=1$ and $\lB(\mb a)$ is stable.
\end{corollary}
\begin{proof}
We know that the sheaf $F_m$ has a filtration with factors
isomorphic to $\lO_C$. Therefore the sheaf $\lB(\mb a,m,\la)$ has
a filtration with factors isomorphic to $\lB(\mb a,1,\la)$. Hence
$\lB(\mb a,m,\la)$ is semistable if and only if $\lB(\mb a,1,\la)$
is semistable and $\lB(\mb a,m,\la)$ can be stable just if $m=1$.
It is clear that $\lB(\mb a,1,\la)$ is (semi)stable if and only if
$\lB(\mb a,1,1)$ is.
\end{proof}

Let us now describe the non-locally free indecomposable sheaves.
Given a chain $\mb a=(\line a1r)$, we define the torsion-free
sheaf $\lS(\mb a)$ as follows. Take the direct image $\pi_*$ of
the sum of $B_i=\lO(a_i)$ and make the following identifications
over the singular point: glue $B_1(p^*)$ with $B_2(p)$, glue
$B_2(p^*)$ with $B_3(p)$ and so on, identifying their bases. The
fibers $B_1(p)$ and $B_r(p^*)$ are not identified.

\begin{theorem}
The sheaves $\lS(\mb a)$ describe all indecomposable torsion-free
non-locally free sheaves over $C$. Different chains induce
non-isomorphic sheaves.
\end{theorem}

Our goal is to determine which of the sheaves $\lB(\mb a,m,\la)$
and $\lS(\mb a)$ are (semi)stable. It follows from Corollary
\ref{crl:lf:sheaves:reduction} that in the case of locally free
sheaves we can restrict ourselves just to $\lB(\mb a)$.

\begin{remark}\label{rmr:degree:rank}
Given a chain $\mb a=(\line a1r)$, one can easily show that $\deg
\lB(\mb a)=\sum a_i$, $\deg\lS(\mb a)=\sum a_i+1$ and $\rank
\lB(\mb a)=\rank\lS(\mb a)=r$.
\end{remark}

\begin{proposition}\label{prp:cyc:exactseq}
Let $\mb a=(\line a1r)$ be a cycle and $\mb b=(\line b1k)$ be its
subchain (see Definition \ref{def:subchain}). Then there is an
exact sequence
$$0\ar\lS((b_1-1,\line b2{k-1},b_k-1))\ar \lB(\mb a)\ar \lS(\mb{b'})\ar 0,$$
where for $k=1$ we consider just $\lS((b_1-2))$ and $\mb{b'}$ is
the complement of $\mb b$ in $\mb a$ .
\end{proposition}
\begin{proof}
Without loss of generality we may assume $\mb b=(\line a1k)$. Let
us denote $B_i=\lO_{\wtl C}(a_i)$. We consider the direct image
$\pi_*$ of the sum
$$(B_1\ts \lO_{\wtl C}(-p))\oplus B_2\oplus\dots\oplus B_{k-1}
\oplus (B_k\ts \lO_{\wtl C}(-p^*))$$ and identify their fibers
precisely like in the construction of \lS. The module obtained in
this way is isomorphic to $\lS((b_1-1,\line b2{k-1},b_k-1))$ and
there is a natural embedding of this module to $\lB(\mb a)$ (the
fiber of $(B_1\ts \lO_{\wtl C}(-p))$ in point $p$ goes to zero
both in fibers $B_r(p^*)$ and $B_1(p)$, and analogously for the
fiber of $(B_k\ts \lO_{\wtl C}(-p^*))$ in point $p^*$). It is
clear that the quotient is isomorphic to the direct image of
$B_{k+1}\oplus\dots\oplus B_r$ with identifications
$B_{k+1}(p^*)\iso B_{k+2}(p),\dots B_{r-1}(p^*)\iso B_r(p)$. But
such a module is precisely $\lS((\line a{k+1}r))$.
\end{proof}

\begin{corollary}\label{crl:cyc:onedirect}
Let $\mb a=(\line a1r)$ be a cycle such that $\lB(\mb a)$ is a
semistable sheaf. Then for any proper subchain $\mb b=(\line b1k)$
of \mb a it holds
$$\frac{\sum_{i=1}^kb_i-1}{k}\le\frac{\sum_{i=1}^ra_i}{r}.$$
If $\lB(\mb a)$ is stable then the inequalities are strict.
\end{corollary}

\begin{proposition}\label{prp:ch:exactseq}
Let $\mb a=(\line a1r)$ be a chain and $\mb b=(\line b1k)$ be its
subchain that does not contain $a_1$ and $a_r$. Then there is an
embedding
$$\lS((b_1-1,\line b2{k-1},b_k-1))\emb \lS(\mb a).$$
If \mb b is a subchain containing $a_1$ or $a_r$ (say, $\mb
b=(\line a1k)$) then there is an exact sequence
$$0\ar\lS((a_1,\line a2{k-1},a_k-1))\ar \lS(\mb a)
\ar \lS((\line a{k+1}r))\ar0.$$
\end{proposition}
\begin{proof}
The proof goes through the same lines as the proof of Proposition
\ref{prp:cyc:exactseq}.
\end{proof}

\begin{corollary}\label{crl:ch:onedirect}
Let $\mb a=(\line a1r)$ be a chain such that $\lS(\mb a)$ is a
semistable sheaf and let
$$\mb{a'}:=(a_1+1,\line a2{r-1},a_r+1).$$
Then for any proper subchain $\mb b=(\line b1k)$ of $\mb {a'}$ it
holds
$$\frac{\sum_{i=1}^kb_i-1}{k}\le\frac{\sum_{i=1}^ra'_i-1}{r}.$$
If $\lS(\mb a)$ is stable then the inequalities are strict.
\end{corollary}
\begin{proof}
If $\lS(\mb a)$ is semistable then for any subchain $\mb b=(\line
b1k)$ of $\mb {a'}$ that does not contain $a'_1$ and $a'_r$ we
have (see Proposition \ref{prp:ch:exactseq} and Remark
\ref{rmr:degree:rank})
$$\frac{\sum_{i=1}^kb_i-2+1}{k}\le\frac{\sum_{i=1}^ra_i+1}{r}=\frac{\sum_{i=1}^ra'_i-1}{r}.$$
If \mb b is a subchain of $\mb{a'}$ containing $a'_1$ or $a'_r$
(say, $\mb b=(\line {a'}1k)$) then according to Proposition
\ref{prp:ch:exactseq} there is an embedding $\lS((b_1-1,\line
b2{k-1},b_k-1))\emb \lS(\mb a)$ and this implies
$$\frac{\sum_{i=1}^kb_i-2+1}{k}\le\frac{\sum_{i=1}^ra_i+1}{r}=\frac{\sum_{i=1}^ra'_i-1}{r}.$$
The claim about stability is analogous.
\end{proof}

Corollaries \ref{crl:cyc:onedirect} and \ref{crl:ch:onedirect}
suggest that one can define stability conditions directly for
chains and cycles. We will do this in the next section. After the
classification of (semi)stable chains and cycles we will be able
to prove that the stability of chains and cycles is not only
necessary for the stability of the corresponding sheaves (as it is
proved in Corollaries \ref{crl:cyc:onedirect} and
\ref{crl:ch:onedirect}), but is also sufficient.

\section{Semistable chains and cycles}\label{sec:ss:ch:cyc}
Recall from Definition \ref{def:ch:cyc} that chains are finite
sequences of integers and cycles are equivalence classes of chains
with respect to the cyclic shift.

\begin{definition}\label{def:subchain}
Given a chain $(\line a1r)$, define its subchain as any chain of
the form $(a_i,a_{i+1},\dots,a_j)$, where $1\le i\le j\le r$.
Given a cycle $(\line a1r)$, define its subchain as any chain of
the form $(a_i,a_{i+1},\dots,a_{i+k})$, where $1\le i\le r$, $0\le
k<r$ and we identify $a_{r+1}$ with $a_1$, $a_{r+2}$ with $a_2$
and so on.
\end{definition}

For example, the cycle $(1,2,3,1,2,3)$ contains the subchain
$(3,1,2,3,1)$ but the chain $(1,2,3,1,2,3)$ does not.

\begin{definition}
Given a chain $\mb a=(\line a1r)$, we call any of its subchains
containing $a_1$ or $a_r$ an extreme subchain of $\mb a$.
\end{definition}

\begin{definition}\label{def:aperiodic}
A cycle is called aperiodic if its sequence cannot be written as a
concatenation of equal proper subsequences.
\end{definition}

%For example, the cycle $(1,2,3,1,2,3)$ is periodic and the cycle $(1,2,3,1,2,3,1,2)$ is aperiodic.

\begin{definition}
Given a chain $\mb a=(\line a1r)$, define its degree, rank, and
slope by
$$\deg\mb a=\sum_{i=1}^r a_i -1,\qquad \rank\mb a=r,\qquad
\mu(\mb a)=\frac{\deg\mb a}{\rank\mb a}.$$ Given a cycle $\mb
a=(\line a1r)$, define its degree, rank, and slope by
$$\deg\mb a=\sum_{i=1}^r a_i ,\qquad \rank\mb a=r,\qquad
\mu(\mb a)=\frac{\deg\mb a}{\rank\mb a}.$$
\end{definition}

For example, the slope of the chain $(1,2,3,1,2,3)$ equals $\frac
{11}{6}$ and the slope of the cycle $(1,2,3,1,2,3)$ equals $2$.

\begin{definition}\label{def:ss:ch:cyc}
The chain (cycle) $\mb a=(\line a1r)$ is called semistable if for
any its subchain \mb b it holds
$$\mu(\mb b)\le \mu(\mb a).$$
If the inequality is strict for any proper subchain then \mb a is
called stable. A proper subchain \mb b of a chain (cycle) \mb a is
called a destabilizing subchain of \mb a if $\mu(\mb b)\ge\mu(\mb
a)$.
\end{definition}

For example, the chain $(1,0,0,1)$ is stable and the cycle
$(1,0,0,1)$ is not stable, because it has slope $1/2$ and contains
a destabilizing subchain $(1,1)$ having the same slope. In what
follows, we will classify (semi)stable chains and cycles. For
example, the only stable chain of rank $7$ and degree $4$ is
$(1,1,0,1,0,1,1)$ and the only stable (aperiodic) cycle of rank
$7$ and degree $4$ is $(1,0,1,0,1,0,1)$.

A chain (cycle) of rank $r$ and degree $d$ will be said to be of
type $(r,d)$. We want to classify all (semi)stable chains and
aperiodic cycles of a fixed type $(r,d)$. The set of semistable
(respectively, stable) chains of type $(r,d)$ will be denoted by
$\mch^{ss}(r,d)$ (respectively, $\mch^s(r,d)$). The set of
aperiodic semistable (stable) cycles of type $(r,d)$ will be
denoted by $\mcyc^{ss}(r,d)$ ($\mcyc^s(r,d)$). The study of
semistable chains and semistable cycles is quite analogous but we
will deal with them separately.

\begin{lemma}\label{lmm:ch:shift}
A chain $\mb a=(a_1,a_2,\dots a_r)$ is (semi)stable if and only if
the chain $(a_1+1,a_2+1\dots,a_r+1)$ is (semi)stable. In
particular, there is a natural bijection $\mch^{ss}(r,d)\iso
\mch^{ss}(r,d+r)$ and we can always assume $0\le\deg\mb a<r$.
\end{lemma}

\begin{lemma}\label{lmm:ch:coker:ineq}
For any subchain $\mb b$ of a semistable chain $\mb
a=(a_1,\dots,a_r)$ one has
$$\mu(\mb a)\le \frac{\sum b_i+1}{\rk\mb b}.$$
Moreover, if $\mb b$ is an extreme subchain then
$$\mu(\mb a)\le \frac{\sum b_i}{\rk\mb b}.$$
\end{lemma}
\begin{proof}
Let us first prove the assertion for extreme subchains. We may
assume that $\mb b=(a_1,\dots,a_{k_1})$. Denote $k_2=r-k_1$,
$x_1=\sum_{i=1}^{k_1}a_i$, $x_2=\sum_{i=k_1+1}^{r}a_i$. Then the
semistability of $\mb a$ implies
$$\frac{x_2-1}{k_2}\le\frac{x_1+x_2-1}{k_1+k_2},$$
hence
$$\frac{x_1+x_2-1}{k_1+k_2}\le\frac{x_1}{k_1}.$$
Let us now assume that $\mb b=(a_{k_1+1},\dots,a_{k_1+k_2})$ is
not extreme, i.e., $k_1\ge1$ and $k_3:=r-k_1-k_2\ge1$. We denote
$x_1=\sum_{i=1}^{k_1}a_i$, $x_2=\sum_{i=k_1+1}^{k_1+k_2}a_i$, and
$x_3=\sum_{i=k_1+k_2+1}^{r}a_i$. It holds by our assumptions
$$\frac{x_1-1}{k_1}\le\frac{x_1+x_2+x_3-1}{k_1+k_2+k_3},\quad
\frac{x_3-1}{k_3}\le\frac{x_1+x_2+x_3-1}{k_1+k_2+k_3},$$ hence
$$\frac{x_1+x_3-2}{k_1+k_3}\le\frac{x_1+x_2+x_3-1}{k_1+k_2+k_3}$$
and therefore
$$\frac{x_1+x_2+x_3-1}{k_1+k_2+k_3}\le\frac{x_2+1}{k_2}.$$
\end{proof}

\begin{corollary}\label{crl:ch:coker:ineq}
If a chain $\mb a=(\line a1r)$ is semistable then $\mu(\mb a)\le
a_1$, $\mu(\mb a)\le a_r$ and for any $k$ one has $\mu(\mb a)\le
a_k+1$.
\end{corollary}

It follows that in a semistable chain $(\line a1r)$ for any
indices $i,j$ one has $a_i-1\le\mu(a)\le a_j+1$ and therefore the
difference between any $a_i$ and $a_j$ is not greater than $2$.
Hence, the elements of \mb a can take at most $3$ consecutive
values.

\begin{lemma}\label{lmm:ch:diff2}
If a semistable chain $\mb a=(\line a1r)$ of type $(r,d)$ contains
elements with difference $2$ then $d$ is a multiple of $r$ (hence
$d=0$ under the assumption $0\le d<r$).
\end{lemma}
\begin{proof}
Assume there are elements in \mb a equal to $m-1$ and $m+1$. Then
we have $(m+1)-1\le\mu(\mb a)\le (m-1)+1$ and therefore $d/r=m$ is
an integer.
\end{proof}

We prove now the first part of Theorem \ref{intro:main:thr2} for
chains. It serves as a basis of our reduction of chains.

\begin{proposition}\label{lmm:ch:reduction}
Let $0<d<r$. Then there is a natural bijection between
$\mch^{ss}(r,d)$ and $\mch^{ss}(d,d-r)$. Analogous with stable
chains.
\end{proposition}
\begin{proof}
Let $\mb a=(\line a1r)$ be a semistable chain of type $(r,d)$. It
follows from Lemma \ref{lmm:ch:diff2} that its elements can take
at most two consecutive values. Obviously, they can be only $0$ or
$1$. From the inequality $a_1\ge\mu(\mb a)>0$ one gets $a_1=1$.
Analogously $a_r=1$. From the condition $\sum_{i=1}^ra_i-1=d$ we
obtain that there are $d+1$ $1$'s among the elements of \mb a. Let
$\line b1d$ be the lengths of consecutive zero-blocks between the
$1$'s. We have $\sum_{i=1}^db_i=r-d-1$. Now, the chain \mb a
consisting of $0$'s and $1$'s is semistable if and only if the
inequality from the definition \ref{def:ss:ch:cyc} holds for any
subchain starting and ending with a one. This can be written as
follows. Any subchain $(\line b{j+1}{j+k})$ of the chain $(\line
b1d)$ should satisfy
$$\frac{(k+1)-1}{\sum_{i=j+1}^{j+k}b_i+k+1}\le\frac{d}{\sum_{i=1}^db_i+d+1},$$
or, equivalently,
$$\frac{\sum_{i=j+1}^{j+k}b_i+1}{k}\ge\frac{\sum_{i=1}^db_i+1}{d},$$
which can be written in the form
$$\frac{\sum_{i=j+1}^{j+k}(-b_i)-1}{k}\le\frac{\sum_{i=1}^d(-b_i)-1}{d}.$$
But this says precisely that the chain $(-b_1,-b_2,\dots,-b_d)$ is
semistable. Its degree is $\sum_{i=1}^d(-b_i)-1=-(r-d-1)-1=d-r$.
The last thing to prove is that, conversely, any such semistable
chain will give nonnegative numbers $b_i$ so that we can
reconstruct the chain \mb a. But the semistability condition for
$(-b_1,\dots,-b_d)$ implies $-b_i-1\le(d-r)/{d}<0$ and therefore
$b_i\ge 0$. This altogether implies that there is a bijection
between $\mch^{ss}(r,d)$ and $\mch^{ss}(d,d-r)$. The proof for
stable chains goes through the same lines.
\end{proof}

This proposition shows that we can reduce the classification of
(semi)stable chains of type $(r,d)$ to the classification of
(semi)stable chains of type $(d,d-r)$, i.e., of those with a
smaller rank. The latter can be reduced to $\mch^{ss}(d,r_0)$
(respectively, to $\mch^s(d,r_0)$), where $0\le r_0<d$ by Lemma
\ref{lmm:ch:shift}. Repeating these reductions we will finally end
up with $\mch^{ss}(h,0)$ (respectively, $\mch^s(h,0)$), where
$h=\gcd(r,d)$. So, the second part of Theorem
\ref{intro:main:thr2} for chains should be proved (and
classification should be done) only for the type $(h,0)$.

For example, let us describe $\mch^{ss}(7,4)$. We write our
reductions as follows
$$\mch^{ss}(7,4)\iso \mch^{ss}(4,4-7)\iso \mch^{ss}(4,1)\iso \mch^{ss}(1,1-4)\iso \mch^{ss}(1,0).$$
Thus, we take the unique element $(1)\in \mch^{ss}(1,0)$ and
reconstruct the element from $\mch^{ss}(7,4)$ going from the right
to the left in our sequence of isomorphisms. We get $(-2)\in
\mch^{ss}(1,-3)$ and therefore the element of $\mch^{ss}(4,1)$
consists of two ones with a zero-block of length $2$ between them,
so we get $(1,0,0,1)\in \mch^{ss}(4,1)$. Then $(0,-1,-1,0)\in
\mch^{ss}(4,-3)$ and the element of $\mch^{ss}(7,4)$ consists of
five ones with zero-blocks of lengths $(0,1,1,0)$ between them, so
we get $(1,1,0,1,0,1,1)\in \mch^{ss}(7,4)$.

\begin{proposition}\label{prp:sschains:r0}
The semistable chains of type $(r,0)$ are of the form
$$(0,\dots,0,1,0,\dots,0,-1,0,\dots,0,1,\dots\dots,-1,0,\dots,0,1,0,\dots,0),$$
where $1$ and $-1$ alternate and the zero-blocks are of arbitrary
lengths (the whole sequence must be, of course, of length $r$). If
$r>1$, none of these chains is stable. If $r=1$ there is precisely
one semistable chain $(1)$ and it is stable.
\end{proposition}
\begin{proof}
Let $\mb a=(\line a1r)$ be semistable of type $(r,0)$. Then we
have $a_i-1\le\mu(\mb a)=0$, $a_i+1\ge\mu(\mb a)=0$ and therefore
$-1\le a_i\le 1$. Let there be $k$ elements in \mb a which are
equal to $1$ and $l$ elements which are equal to $-1$. We have
then $\deg\mb a=l-k-1=0$ and so $l=k+1$. If there exists a
subchain \mb b containing only zeros and ones with at least two
ones then $\deg\mb b\ge 2-1>0$ and therefore $\mu(\mb b)>\mu(\mb
a)=0$, which is impossible. This together with $l=k+1$ imply that
$1$ and $-1$ alternate in \mb a and therefore \mb a has a required
form. Conversely, if a chain \mb a has the form like in the
condition of the proposition then, first of all, its degree equals
$0$. For any subchain \mb b the difference between the numbers of
$1$'s and $-1$'s is not greater than $1$ and therefore $\deg\mb
b\le 0$, which implies $\mu(\mb b)\le\mu(\mb a)=0$. To prove that
\mb a is not stable if $r>1$ we notice that for a proper subchain
$(1)$ of \mb a one has $\mu((1))=0=\mu(\mb a)$. The last assertion
of the proposition is trivial.
\end{proof}

This proposition together with Proposition \ref{lmm:ch:reduction}
implies Theorem \ref{intro:main:thr2} for chains. The further
considerations are of independent interest.

\begin{lemma}
A chain $\mb a=(a_1,\dots,a_r)$ is semistable (stable) if and only
if for any of its extreme subchains \mb b it holds $\mu(\mb
b)\le\mu(\mb a)$ ($\mu(\mb b)<\mu(\mb a)$). In particular, if a
chain \mb a is non-stable then it contains an extreme
destabilizing subchain.
\end{lemma}
\begin{proof}
Assuming that for any extreme subchain \mb b of  \mb a it holds
$\mu(\mb b)\le\mu(\mb a)$, we will show that $\mb a$ is
semistable. Let $\mb c=(a_{k_1+1},a_{k_1+2},\dots a_{k_1+k_2})$ be
a subchain of \mb a. We denote $k_3=r-k_1-k_2$,
$x_1=\sum_{i=1}^{k_1}a_i$, $x_2=\sum_{i=k_1+1}^{k_1+k_2}a_i$, and
$x_3=\sum_{i=k_1+k_2+1}^{r}a_i$. We want to show that $\mu(\mb
c)\le\mu(\mb a)$, so we may suppose that \mb c is not extreme,
hence $k_1\ge 1$ and $k_3\ge 1$.
%It holds by our assumption
%$$\frac{x_1+x_2-1}{k_1+k_2}\le\frac{x_1+x_2+x_3-1}{k_1+k_2+k_3},\quad
%\frac{x_2+x_3-1}{k_2+k_3}\le\frac{x_1+x_2+x_3-1}{k_1+k_2+k_3},$$ which imply
The same proof as in Lemma \ref{lmm:ch:coker:ineq} shows
$$\frac{x_1+x_2+x_3-1}{k_1+k_2+k_3}\le\frac {x_3}{k_3},\quad
\frac{x_1+x_2+x_3-1}{k_1+k_2+k_3}\le\frac {x_1}{k_1}$$ and
therefore
$$\frac{x_1+x_2+x_3-1}{k_1+k_2+k_3}\le\frac{x_1+x_3}{k_1+k_3}.$$
This implies
$$\frac{x_2-1}{k_2}\le\frac{x_1+x_2+x_3-1}{k_1+k_2+k_3},$$
i.e., $\mu(\mb c)\le\mu(\mb a)$. The proof for stable chains is
analogous.
\end{proof}

\begin{lemma}\label{lmm:ch:factors}
Let $\mb a=(\line a1r)$ be a semistable chain and $\mb b=(\line
a1k)$ be its extreme destabilizing subchain. Then the chains $\mb
b$ and $\mb{b'}=(a_{k+1}+1,\line a{k+2}r)$ are semistable chains
with slope $\mu(\mb a)$.
\end{lemma}
\begin{proof}
It follows from the condition
$$\frac{\sum_{i=1}^ka_i-1}k=\frac{(\sum_{i=1}^ka_i-1)+\sum_{i=k+1}^ra_i}r$$
that
$$\mu(\mb{b'})=\frac{\sum_{i=k+1}^ra_i}{r-k}=
\frac{(\sum_{i=1}^ka_i-1)+\sum_{i=k+1}^ra_i}r=\mu(\mb a).$$ The
semistability of \mb b is trivial. To prove the semistability of
$\mb{b'}$ we note that if \mb c is a subchain of $\mb b'$ not
containing the element $a_{k+1}+1$ then $\mu(\mb c)\le\mu(\mb
a)=\mu(\mb{b'})$. If \mb c contains $a_{k+1}+1$ then it is of the
form $(a_{k+1}+1,\line a{k+2}{k+l})$ and therefore it would follow
from
$$\mu(\mb c)=\frac{\sum_{i=k+1}^{k+l}a_i}{l}>\mu(\mb a),\qquad
\frac{\sum_{i=1}^ka_i-1}{k}=\mu(\mb a)$$ that
$$\frac{\sum_{i=1}^{k+l}a_i-1}{k+l}>\mu(\mb a),$$
which is impossible as \mb a is semistable.
\end{proof}

We return to (semi)stable cycles.

\begin{lemma}
The cycle $\mb a=(a_1,a_2,\dots a_r)$ is semistable (stable) if
and only if the cycle $(a_1+1,a_2+1\dots,a_r+1)$ is semistable
(stable). In particular, there is a bijection $\mcyc^{ss}(r,d)\iso
\mcyc^{ss}(r,d+r)$ and we may always assume $0\le\deg a<r$.
\end{lemma}

\begin{lemma}
If the chain $\mb a=(\line a1r)$ is semistable then for any index
$k$ it holds $\mu(\mb a)\le a_k+1$.
\end{lemma}
\begin{proof}
Without loss of generality we may assume $k=r$. Semistability of
\mb a implies
$$(\sum_{i=1}^{r-1}a_i-1)/(r-1)\le(\sum_{i=1}^{r}a_i)/r,$$
hence
$$\sum_{i=1}^ra_i\le ra_r+r$$
and the claim follows.
\end{proof}

It follows that for any semistable cycle $\mb a=(\line a1r)$ and
any indices $i,j$ it holds $a_i-1\le\mu(\mb a)\le a_j-1$. As
above, we obtain that the elements of \mb a can take at most three
consecutive values.

\begin{lemma}
If a semistable cycle of type $(r,d)$ contains elements with
difference $2$ then $d$ is a multiple of $r$ (hence $d=0$ under
the assumption $0\le d<r$).
\end{lemma}
\begin{proof}
The proof is the same as the proof of Lemma \ref{lmm:ch:diff2}
\end{proof}

\begin{proposition}\label{prp:cyc:reduction}
Let $0<d<r$. Then there is a bijection between $\mcyc^{ss}(r,d)$
and $\mcyc^{ss}(d,d-r)$. Analogous with stable aperiodic cycles.
\end{proposition}
\begin{proof}
Let $\mb a=(\line a1r)$ be a semistable cycle of type $(r,d)$. We
know that its elements can take at most two consecutive values.
Obviously, they can be only $0$ and $1$. From the condition
$\sum_{i=1}^ra_i=d$ we get that there are $d$ ones among the
elements of \mb a. Let $\line b1d$ be the lengths of consecutive
zero-blocks between the ones. We have $\sum_{i=1}^db_i=r-d$. Now,
the cycle \mb a consisting of zeros and ones is semistable if and
only if the inequality from Definition \ref{def:ss:ch:cyc} holds
for any subchain starting and ending with a one. This can be
written as follows. For any subchain $(\line b{j+1}{j+k})$ of the
cycle $(\line b1d)$ one should have
$$\frac{(k+1)-1}{\sum_{i=j+1}^{j+k}b_i+k+1}\le\frac{d}{\sum_{i=1}^db_i+d},$$
or, equivalently,
$$\frac{\sum_{i=j+1}^{j+k}b_i+1}{k}\ge\frac{\sum_{i=1}^db_i}{d},$$
which can be written in the form
$$\frac{\sum_{i=j+1}^{j+k}(-b_i)-1}{k}\le\frac{\sum_{i=1}^d(-b_i)}{d}.$$
But this says precisely that the cycle $(-b_1,-b_2,\dots,-b_d)$ is
semistable. Its degree equals $\sum_{i=1}^d(-b_i)=d-r$. It remains
to prove that, conversely, any such semistable cycle will produce
nonnegative numbers $b_i$ so that we can reconstruct the cycle \mb
a. But the semistability condition implies $-b_i-1\le(d-r)/{d}<0$,
therefore $b_i\ge 0$. It is clear that the cycle \mb a is
aperiodic if and only if \mb b is aperiodic. Altogether it implies
that there is a bijection between $\mcyc^{ss}(r,d)$ and
$\mcyc^{ss}(d,d-r)$. The proof for stable cycles goes through the
same lines.
\end{proof}

Using this proposition, precisely as it was done for chains, we
can reduce the study of $\mcyc^{ss}(r,d)$ to the study of
$\mcyc^{ss}(h,0)$, where $h=\gcd(r,d)$.

For example, let us describe $\mcyc^{ss}(7,4)$. We write our
reductions as follows
$$\mcyc^{ss}(7,4)\iso \mcyc^{ss}(4,4-7)\iso \mcyc^{ss}(4,1)\iso \mcyc^{ss}(1,1-4)\iso \mcyc^{ss}(1,0).$$
Thus, we take the unique element $(0)\in \mcyc^{ss}(1,0)$ and
reconstruct the element from $\mcyc^{ss}(7,4)$ going from the
right to the left in our sequence of isomorphisms. We get $(-3)\in
\mcyc^{ss}(1,-3)$ and therefore the element of $\mcyc^{ss}(4,1)$
equals $(1,0,0,0)$ (the length of zero-block equals $3$). Then
$(0,-1,-1,-1)\in \mcyc^{ss}(4,-3)$ and the element of
$\mcyc^{ss}(7,4)$ has zero-blocks of lengths $(0,1,1,1)$, so it
looks like $(1,1,0,1,0,1,0)$. Clearly, it is equivalent to
$(1,0,1,0,1,0,1)\in \mcyc^{ss}(7,4)$.

\begin{proposition}\label{prp:sscycles:r0}
The semistable cycles of type $(r,0)$ are of the form
$$(0,\dots,0,1,0,\dots,0,-1,0,\dots,0,1,\dots\dots,-1,\dots,0),$$
where $1$ and $-1$ alternate and zero-blocks are arbitrary (the
sequence should of course be of length $r$). If $r>1$, none of
these cycles is stable aperiodic. If $r=1$ there is just one
semistable cycle $(1)$ and it is stable.
\end{proposition}
\begin{proof}
Let $\mb a=(\line a1r)$ be semistable of type $(r,0)$. Then we
have $a_i-1\le\mu(\mb a)=0$, $a_i+1\ge\mu(\mb a)=0$ and therefore
$-1\le a_i\le 1$. Let there be $k$ elements in \mb a which equal
$1$ and $l$ elements which equal $-1$. We have then $\deg\mb
a=l-k=0$, so $l=k$. If there exists a subchain \mb b containing
only zeros and ones with at least two ones then $\deg\mb b\ge
2-1>0$ and therefore $\mu(\mb b)>\mu(\mb a)=0$, which is
impossible. This, with $l=k$ imply that $1$ and $-1$ alternate in
\mb a and therefore \mb a has the required form. Conversely, if a
chain \mb a has the form like in the condition of the proposition
then, first of all, its degree equals $0$. For any subchain \mb b
the difference between the numbers of $1$'s and $-1$'s is no
greater than $1$ and therefore $\deg\mb b\le 0$, which implies
$\mu(\mb b)\le\mu(\mb a)=0$. To prove that any aperiodic \mb a is
non-stable if $r>1$ we notice that it contains nonzero elements,
because otherwise it would be periodic. But for a proper subchain
$(1)$ of \mb a one has $\mu((1))=0=\mu(\mb a)$, so \mb a is
non-stable. The last assertion of the proposition is trivial.
\end{proof}

This proposition together with Proposition \ref{prp:cyc:reduction}
implies Theorem \ref{intro:main:thr2} for cycles.

\begin{lemma}\label{lmm:cyc:factors}
Let $\mb a=(\line a1r)$ be a semistable cycle and $\mb b=(\line
a1k)$ be its destabilizing subchain. Then the chains $\mb b$ and
$\mb{b'}=(a_{k+1}+1,\line a{k+2}{r-1},a_r+1)$ are semistable
chains with slope $\mu(\mb a)$.
\end{lemma}
\begin{proof}
It follows from the condition
$$\frac{\sum_{i=1}^ka_i-1}k=\frac{(\sum_{i=1}^ka_i-1)+(\sum_{i=k+1}^ra_i+1)}r$$
that
$$\mu(\mb{b'})=\frac{\sum_{i=k+1}^ra_i+1}{r-k}=
\frac{(\sum_{i=1}^ka_i-1)+(\sum_{i=k+1}^ra_i+1)}r=\mu(\mb a).$$
The semistability of \mb b is trivial. To prove the semistability
of $\mb{b'}$ we note that if \mb c is a subchain of $\mb b'$ not
containing elements $a_{k+1}+1$ and $a_{r}+1$ then $\mu(\mb
c)\le\mu(\mb a)=\mu(\mb{b'})$. If \mb c is a proper subchain of
$\mb{b'}$ containing, say, $a_{k+1}+1$ then it is of the form
$(a_{k+1}+1,\line a{k+2}{k+l})$ and therefore it would follow from
$$\mu(\mb c)=\frac{\sum_{i=k+1}^{k+l}a_i}{l}>\mu(\mb a),\qquad
\frac{\sum_{i=1}^ka_i-1}{k}=\mu(\mb a)$$ that
$$\frac{\sum_{i=1}^{k+l}a_i-1}{k+l}>\mu(\mb a),$$
which is impossible as \mb a is semistable.
\end{proof}

\section{Classification of semistable sheaves}\label{sec:ss:sh:classif}
We know how to classify the (semi)stable chains and cycles, so the
classification of (semi)stable sheaves will be complete if we will
prove that it holds the converse of Corollaries
\ref{crl:cyc:onedirect} and \ref{crl:ch:onedirect}. We do this in
four steps.

\begin{lemma}
The sheaf $\lB(\mb a)$ is stable if and only if the cycle $\mb a$
is stable. In this case degree and rank are coprime.
\end{lemma}
\begin{proof}
The only if part is already proved. Let \mb a be a stable cycle of
type $(r,d)$. We know that necessarily $r$ and $d$ are coprime and
$\mb a$ is the unique stable cycle of type $(r,d)$.
%The fact that $M_C(r,d)$ is one-dimensional
%implies that there should exist stable sheaves of the
%form $\lB(\mb b,m,\la)$ (the sheaves of the form $\lS(\mb b)$
%can give just a discrete set).
There exist stable locally free sheaves of type $(r,d)$ (see e.g.\
\cite{Bu}). Let $\lB(\mb b,m,\la)$ be any of them. Then $m=1$ and
\mb b is stable of type $(r,d)$, hence $\mb b=\mb a$. It follows
that $\lB(\mb a)$ is stable.
\end{proof}

\begin{lemma}
The sheaf $\lS(\mb a)$ ($\mb a=(\line a1r)$) is stable if and only
if the chain $\mb{a'}=(a_1+1, \line a2{r-1}, a_r+1)$ is stable. In
this case degree and rank are coprime.
\end{lemma}
\begin{proof}
The only if part is already proved. Let $\mb{a'}$ be a stable
chain of type $(r,d)$. Then $r$ and $d$ are coprime and $\mb a'$
is the unique stable chain of type $(r,d)$. Let $M_C(r,d)$ denote
the moduli space of stable sheaves of type $(r,d)$ over $C$. The
subspace of $M_C(r,d)$ consisting of the locally free sheaves
$\lB(\mb b,1,\la)$ (where \mb b is a unique stable cycle of type
$(r,d)$) is isomorphic to $k^*$. It follows from the projectivity
of $M_C(r,d)$ that it cannot coincide with $k^*$ and therefore it
contains some $\lS(\mb c)$, so that the corresponding chain
$\mb{c'}$ of type $(r,d)$ is stable and we deduce from the
uniqueness of stable chains of type $(r,d)$ that $\mb{c'}=\mb{a'}$
hence $\mb c=\mb a$ and $\lS(\mb a)$ is stable.
\end{proof}

\begin{lemma}
The sheaf $\lS(\mb a)$ ($\mb a=(\line a1r)$) is semistable if and
only if the chain $\mb{a'}=(a_1+1,\line a2{r-1}, a_r+1)$ is
semistable.
\end{lemma}
\begin{proof}
The only if part is already proved. Conversely, if the chain $\mb
a'$ is stable, then we are done. So, let us assume that $\mb a'$
is semistable but not stable. Then it contains an extreme
destabilizing subchain, which without loss of generality we will
assume to be of the form $(a_1+1,\line a2k)$. By Lemma
\ref{lmm:ch:factors}, we know that the chains $(a_1+1,\line a2k)$
and $(a_{k+1}+1,\line a{k+2}{r-1},a_r+1)$ are semistable with the
same slope $\mu(\mb a')$, so by induction on rank we deduce that
$\lS((a_1,a_2,\dots,a_{k-1},a_k-1))$ and
$\lS((a_{k+1},a_{k+2},\dots,a_{r-1},a_r))$ are semistable with the
slope $\mu(\mb {a'})$. Now, it follows from the exact sequence of
Proposition \ref{prp:ch:exactseq} that $\lS(\mb a)$ is also
semistable.
\end{proof}

\begin{lemma}
The sheaf $\lB(\mb a)$ is semistable if and only if the cycle \mb
a is semistable.
\end{lemma}
\begin{proof}
The only if part is already proved. Conversely, if the cycle \mb a
is stable, then we are done. So, let us assume that \mb a is
semistable but not stable. Then it contains a destabilizing
subchain which, without loss of generality, we will assume to be
of the form $(\line a1k)$. By Lemma \ref{lmm:cyc:factors}, we know
that the chains $(a_1,\line a2k)$ and $(a_{k+1}+1,\line
a{k+2}{r-1},a_r+1)$ are semistable with the same slope $\mu(\mb
a)$, therefore the sheaves $\lS((a_1-1,\line a2{k-1},a_k-1))$ and
$\lS((a_{k+1},\line a{k+2}{r-1},a_r))$ are semistable with the
slope $\mu(\mb{a'})$. Now, it follows from the exact sequence of
Proposition \ref{prp:cyc:exactseq} that $\lB(\mb a)$ is also
semistable.
\end{proof}

Altogether it proves Theorem \ref{intro:main:thr1}. We formulate
now some corollaries.

%\newtheorem*{copy-main-proposition}{Proposition \ref{intro:main:thr1}}
%\begin{copy-main-proposition}
%Given an aperiodic cycle \mb a, the sheaf  $\lB(\mb a)$ is
%(semi)stable if and only if the cycle $\mb a$ is (semi)stable.
%Given a chain $\mb b=(\line b1r)$, the sheaf $\lS(\mb b)$ is
%(semi)stable if and only if the chain
%$(b_1+1,b_2,\dots,b_{r-1},b_r+1)$ is (semi)stable.
%\end{copy-main-proposition}

\begin{corollary}
If $\gcd(r,d)>1$ then there are no stable sheaves in $\lE(r,d)$.
The number of non-locally free semistable sheaves in $\lE(r,d)$ is
finite and non-zero. The family of semistable locally free sheaves
in $\lE(r,d)$ is parameterized by a finite (non-empty) union of
copies of $k^*$.
\end{corollary}

\begin{corollary}\label{crl:ssh:properties}
If $\gcd(r,d)=1$ then all semistable sheaves in $\lE(r,d)$ are
stable. There is precisely one non-locally free semistable sheaf.
The family of semistable locally free sheaves in $\lE(r,d)$ is
parameterized by $k^*$
\end{corollary}

\begin{definition}\label{hereditary}
We call a semistable sheaf $F$ homogeneous if all the stable
factors of its Jordan-H\"older filtration are isomorphic. The
corresponding isomorphism class is called a basic block of $F$. In
particular a stable sheaf is homogeneous.
\end{definition}

\begin{corollary}
All indecomposable semistable sheaves over $C$ are homogeneous.
\end{corollary}
\begin{proof}
First of all, we note that given two non-isomorphic stable sheaves
$G_1,G_2$ of the same type $(r,d)$, we have $\Ext^1(G_1,G_2)=0$.
This follows immediately from the Serre duality which is
applicable because one of two sheaves $G_1,G_2$ is necessarily
locally free (see Corollary \ref{crl:ssh:properties}). This
implies, that also for any two homogeneous sheaves $G_1,G_2$
having non-isomorphic basic blocks of the same type one has
$\Ext^1(G_1,G_2)$. Consider a Jordan-H\"older filtration of a
given semistable sheaf $F$. If we have two consecutive
non-isomorphic factors $G_1,G_2$ then we can change the filtration
in such a way that $G_1$ and $G_2$ are interchanged (using
$\Ext^1(G_1,G_2)=0)$. This shows that $F$ has a filtration with
homogeneous factors having pairwise different basic blocks. As we
have shown, the $\Ext^1$-group between two different factors is
zero and therefore such a filtration necessarily splits.
\end{proof}

%GATHER{../tex/fullbib.bib}   % For Gather Purpose Only
\bibliography{fullbib}

\providecommand{\bysame}{\leavevmode\hbox to3em{\hrulefill}\thinspace}
\providecommand{\href}[2]{#2}
\begin{thebibliography}{1}

\bibitem{Bu}
Igor Burban, \emph{Stable bundles on a ratinal curve with one simple double
  point}, Ukr. Math. J. \textbf{55} (2003), 1043--1053.

\bibitem{BK2}
Igor Burban and Bernd Kreu{\ss}ler, \emph{Simple torsion free sheaves on a
  nodal weierstra{\ss} curve}, in preparation.

\bibitem{BK1}
\bysame, \emph{{Fourier-Mukai transforms and semi-stable sheaves on nodal
  Weierstra{\ss} cubics}}, J. Reine Angew. Math. \textbf{584} (2005), 45--82.

\bibitem{DG}
Yu. Drozd and G.-M. Greuel, \emph{Tame and wild projective curves and
  classification of vector bundles}, J. Algebra \textbf{246} (2001), no.~1,
  1--54.

\bibitem{FMW}
R.~Friedman, J.~Morgan, and E.~Witten, \emph{Vector bundles over elliptic
  fibrations}, J. Algebr. Geom. \textbf{8} (1999), 279--401.

\bibitem{GP1}
I.M. Gelfand and V.A. Ponomarev, \emph{Indecomposable representations of the
  {L}orentz group}, Russ. Math. Surv. \textbf{23} (1968), no.~2, 1--58.

\bibitem{Yu}
Ivan Yudin, \emph{Tensor product of vector bundles on configurations of
  projective lines}, 2000, Diploma Thesis, Kaiserslautern.

\end{thebibliography}
\bibliographystyle{hamsplain}
\end{document}